\documentclass[a4paper,11pt]{article}

\usepackage{latexsym}
\usepackage{amsmath,amsfonts,amssymb,amsthm}
\newtheorem{thm}{Theorem}[section]
\newtheorem{lem}[thm]{Lemma}

\newtheorem{prop}[thm]{Proposition}
\newtheorem{conj}[thm]{Conjecture}
\newtheorem{rem}[thm]{Remark}

\title{Counting packings of generic subsets in finite groups}
\author{Roland Bacher}

\begin{document}
\maketitle

{\sl Abstract\footnote{
Math. class: 05A15, 05C30, 11B73, 11P99}: 
A packing of subsets $\mathcal S_1,\dots,
\mathcal S_n$ in a group $G$ is an element $(g_1,\dots,g_n)$ of $G^n$
such that $g_1\mathcal S_1,\dots,g_n\mathcal S_n$ are 
disjoint subsets of $G$. We give a formula for the 
number of packings if the group $G$ is finite and if the subsets
$\mathcal S_1,\dots,\mathcal S_n$ 
satisfy a genericity condition. This formula can be seen as a 
generalization of the falling factorials which encode the number of packings 
in the case where all the sets $\mathcal S_i$ are singletons.}

\section{Introduction}

A {\it (left-)packing} of $n$ non-empty subsets 
$\mathcal S_1,\dots,\mathcal S_n$ 
in a group $G$ is an element $(g_1,\dots,g_n)$ of $G^n$
such that the left-translates 
$g_1\mathcal S_1,\dots,g_n\mathcal S_n$
of the sets $\mathcal S_i$ 
are disjoint. 
The sets $\mathcal S_1,\dots,\mathcal S_n$
are labelled by their indices. In particular, permuting 
the elements $g_1,\dots,g_n$ of a packing $(g_1,\dots,g_n)\in\mathcal G^n$
of $\mathcal S_1=\dots=\mathcal S_n$
yields a different packing. Moreover, in the case where 
$\mathcal S_1$ for example is of the form $\mathcal S_1=H\mathcal S_1$ 
for some subgroup $H$ of $G$, a packing $(g_1,\dots,g_n)$ gives rise to 
$\sharp(H)$ distinct packings $(g_1h,g_2,\dots,g_n),\ h\in H$.

There is an obvious one-to-one map between 
packings of $\mathcal S_1,\dots,\mathcal S_n\subset G$ and packings 
of $a_1\mathcal S_1,\dots,a_n\mathcal S_n\subset G$ for every 
$(a_1,\dots,a_n)\in G^n$.

This paper deals with enumerative properties of left-packings
in the case where $G$ is a finite group.
Using the involutive antiautomorphism $g\longmapsto g^{-1}$,
its content can easily be modified in order to deal with 
right-packings $\mathcal S_1 g_1,\dots,\mathcal S_ng_n$.

In the sequel, we denote by $\alpha(G;\mathcal S_1,\dots,\mathcal S_n)
\leq N^n$
the number of packings of $n$ non-empty subsets $\mathcal S_1,\dots,
\mathcal S_n$ in a finite group $G$ with $N$ elements. 
Computing $\alpha(G;\mathcal S_1,\dots,\mathcal S_n)$ for arbitrary subsets $\mathcal S_1,\dots,\mathcal S_n$
in a finite group $G$ is probably difficult. There are however
easy lower and upper bounds:

\begin{prop}\label{propalphaineq} We set
$a=\alpha(G;\mathcal S_1,\dots,\mathcal S_n)$ and 
$b=\alpha(G;\mathcal S_1,\dots,\mathcal S_n,\mathcal S_{n+1})$
where $\mathcal S_1,\dots,\mathcal S_{n+1}$ are $(n+1)$ non-empty 
subsets in a finite group $G$. We have the inequalities
$$\left(N-\sharp(\mathcal S_{n+1})
\sum_{i=1}^n\sharp(\mathcal S_i)\right)a\leq b\leq
\left(N-\sum_{i=1}^n\sharp(\mathcal S_i)\right)a\ .$$

In particular, we have 
\begin{eqnarray}\label{formeqalpha}
b=\left(N-\sum_{i=1}^n\sharp(\mathcal S_i)\right)a
\end{eqnarray}
if $\mathcal S_{n+1}$ is a singleton.
\end{prop}

Proposition \ref{propalphaineq} will be proven in Section 
\ref{sectgeneralcase}.

A family $\mathcal S_1,\dots,\mathcal S_n$ 
of $n$ non-empty subsets in a group $G$ with identity element $e$
is {\it generic} if for every sequence $i_1,\dots,i_k$ of $k$
distinct elements in $\{1,\dots,n\}$ and for every choice of 
elements $g_{i_j}\in\mathcal S_{i_j}^{-1}\mathcal S_{i_j}\setminus \{e\}$, 
we have
$$g_{i_1}g_{i_2}\cdots g_{i_k}\not= e\ .$$

Otherwise stated, a family $\mathcal S_1,\dots,\mathcal S_n$ of
subsets  in a group $G$ is generic if the only solution of the equations 
$g_{i_1}\cdots g_{i_n}=e$
with $g_{i_j}\in\mathcal S_{i_j}^{-1}\mathcal S_{i_j}$
for $\{i_1,\dots,i_n\}=\{1,\dots,n\}$  is given by 
$g_{i_j}=e$ for all $j$.

Genericity excludes \lq\lq accidental intersections\rq\rq among 
translates $g_1\mathcal S_1,\dots,g_n\mathcal S_n$ in the following sense: 
Given a collection of translates $g_1\mathcal S_1,\dots,g_n\mathcal S_n$,
we consider the associated {\it intersection graph} with vertices
$\mathcal S_i$ and edges joining $\mathcal S_i,\mathcal S_j$ 
if $g_i\mathcal S_i\cap g_j\mathcal S_j\not=\emptyset$.
Genericity of a family $\mathcal S_1,\dots,\mathcal S_j$ in a group 
$G$ is equivalent to the statement that all intersection graphs
are primal graphs of hyperforests. Intuitively speaking, 
intersections among translates
of a generic family are always \lq\lq as small as possible\rq\rq.

\noindent{\bf Example.}
Genericity in an additive abelian group $G$ boils down to the fact that the
subset $(\mathcal S_1-\mathcal S_1)
\times\dots\times(\mathcal S_n-\mathcal S_n)$ of the group $G^n$
intersects the subgroup $\{(x_1,\dots,x_n)\in G^n\ \vert\ 
\sum_{i=1}^nx_i=0\}$ of $G^n$ only in the 
identity element $(0,\dots,0)$.

A generic family $\mathcal S_1,\dots,\mathcal S_n$ of subsets in 
the additive group $\mathbb Z$ 
with prescribed cardinalities $s_i=\sharp (\mathcal S_i)$ can be 
constructed by starting with 
$\mathcal S_1=\{0,\dots,s_1-1\}$ and by 
defining $\mathcal S_i$ recursively as
$\mathcal S_i=\{0,k_i,2k_i,\dots,(s_i-1)k_i\}$
where $k_i$ is an arbitrary natural integer strictly larger than 
$\sum_{j=1}^{i-1}\left(\max(\mathcal S_j)-\min(\mathcal S_j)\right)=
\sum_{j=1}^{i-1}(s_j-1)k_j$.
A generic family is thus for example given by the sets 
$\mathcal S_1=\{0,1\},\mathcal S_2=\{0,2\},\dots,\mathcal S_i=\{0,2^{i-1}\},
\dots,\mathcal S_n=\{0,2^{n-1}\}$.

Reduction of a generic family $\mathcal S_1,\dots,\mathcal S_n
\subset \mathbb Z$ modulo a natural integer
$N$ yields a generic family of $\mathbb Z/N\mathbb Z$ 
except if $N$ is a divisor of a non-zero integer
in the finite set $\{\sum_{i=1}^n \left(\mathcal S_i-\mathcal S_i
\right)\}$.

\begin{rem}
  The terminology ``generic family'' can be motivated as follows: 
Given $n$ strictly positive natural numbers $s_1,\dots,s_n$, 
most uniform random 
choices of $n$ subsets $\mathcal S_1,\dots,\mathcal S_n$ with 
$\sharp (\mathcal S_i)=s_i$ (among all ${N\choose s_i}$
possible subsets) in a finite group $G$ of order $N$ 
should yield a generic family if $N$ is large compared to 
$\sum_{k=2}^nk!\tau_k$
with $\tau_2,\dots,\tau_n$ defined by $\sum_{k=0}^n \tau_kt^k=\prod_{j=1}^n
\left(1+s_j(s_j-1)t\right)$. Indeed, the number
$k!\tau_k$ is an upper bound on the number 
of elements in the set $E_k$ containing all products of the form 
$g_{i_1}\cdots g_{i_k}$ with $g_{i_j}\in \mathcal S_{i_j}^{-1}\mathcal S_{i_j}
\setminus\{e\}$ and $i_1,\dots,i_k$ given by $k\geq 2$ distinct elements 
of $\{1,\dots,n\}$. Under the (naive but hopefully correct) assumption that the elements of 
$E_k$ are uniformly distributed in $G$, the probability for non-genericity of 
$\mathcal S_1,\dots,\mathcal S_n$ is at most 
$\frac{1}{N}\sum_{k=2}^nk!\kappa_k$.
Observe also the trivial inequalities
$\sum_{k=2}^nk!\kappa_k<n!\sum_{k=0}^n\kappa_k=n!\prod_{j=1}^n(1+s_j(s_j-1))
\leq n!\prod_{j=1}^ns_j^2$.
\end{rem}

The aim of this paper is to describe a universal formula
for the number of packings for a generic family of subsets 
$\mathcal S_1,\dots,\mathcal S_n$
in a finite group $G$.
The number of associated packings depends then only on the 
cardinalities of $G$ and  $\mathcal S_1,\dots,\mathcal S_n$.
Moreover for fixed cardinalities of $\mathcal S_1,\dots,\mathcal S_n$,
the dependency on the cardinality of $G$ is polynomial of degree $n$.
A trivial example is the generic family
given by $n$ subsets reduced to singletons. The associated number of 
packings in a finite group with $N$ elements is then easily seen
to be given by the polynomial $n!{N\choose n}=N(N-1)\cdots (N-n+1)
\in\mathbb Z[N]$
with coefficients given by Stirling 
numbers of the first kind. This polynomial is also called a {\em
falling factorial} and denoted by $N^{\underline n}$. Using the formulae
of our paper, it is possible to define the falling factorial 
$N^{\underline \lambda}$ associated to a partition $\lambda=\lambda_1,\lambda_2,\dots$ by counting packings of generic families with $\lambda_1$ 
subsets having $\nu_1,\nu_2,\dots,\nu_{\lambda_1}$ elements 
where $\nu_i=\{j\ \vert \lambda_j\geq i\}$ is the $i-$th part of the 
transposed partition $\nu=\lambda^t$ of $\lambda$.
The map $\lambda\longmapsto N^{\underline \lambda}$ is however perhaps 
not exceedingly interesting. On one hand, it is not into
since $N^{\underline\lambda}=N$ for every partition $\lambda$ of the form 
$1,1,1,\dots$. On the other hand, fixing the content $\sum_j \lambda_j$
of the partition $\lambda$, our formulae show that the coefficients 
of $N^{\underline\lambda}$ depend linearly on the elementary symmetric functions 
$\sigma_2=\sum_{i<j}\nu_i\nu_j,\sigma_3=\sum_{i<j<k}\nu_i\nu_j\nu_k,\dots,
\sigma_{\lambda_1}=\nu_1\nu_2\cdots \nu_{\lambda_1}$
of the partition $\nu=\lambda^t$.

The study of generic packings in groups is, as far as I am aware,
a new addition to the already large set of classical 
notions of packings. Well-known and well-studied examples are 
lattice-packings in Euclidean spaces or more generally
sphere-packings in metric spaces. 
Error-correcting codes corresponding to packings 
of spheres (with respect to the Hamming distance given by the 
number of distinct coordinates) into $\mathbb F_q^d$ are discrete analogues.
The associated theories
have however a different flavour since one tries 
to pack a huge (perhaps infinite) number of identical copies of
spheres as tightly as possible. 

Subsets in generic families are in general all distinct:
Repetition destroys genericity except in the case of singletons.
Moreover, packings of generic sets have typically very small densities. 
Generic families are mainly interesting 
for enumerative properties of the corresponding packings.

This paper is organized as follows: Section \ref{sectmain}
contains the main result, Theorem \ref{thmU}.
It expresses the number of packings
of a generic family $\mathcal S_1,\dots,\mathcal S_n$ in a finite group
in terms of a formal power series 
$U=U(x,\sigma_1,\sigma_2,\dots)\in A[[x]]$ with coefficients
in the ring $A=\mathbb Z[\sigma_1,\sigma_2,\dots]$
of polynomials in elementary symmetric functions $\sigma_1
=\sum_{i=1}^n \sharp(\mathcal S_i),\sigma_2=
\sum_{i<j}\sharp(\mathcal S_i)\sharp(\mathcal S_j),\dots$
of $\sharp(\mathcal S_1),\dots,\sharp(\mathcal S_n)$.
The series $U$ is given
explicitly by Formula (\ref{formulaU}) and involves combinatorial 
integers $t_{i,j}(n)$ (defined recursively by Formula 
(\ref{eqtriangle})) which extend Stirling numbers of the first kind.
The first few coefficients of $U$ are given by
\begin{eqnarray*}
&&1-\sigma_2 x-((1-\sigma_1)\sigma_3+\sigma_4)x^2\\
&&-((2-3\sigma_1+\sigma_1^2)\sigma_4+(5-3\sigma_1)\sigma_5+3\sigma_6)x^3\\
&&-((6-11\sigma_1+6\sigma_1^2-\sigma_1^3)\sigma_5+(26-26\sigma_1+6\sigma_2^2)
\sigma_6\\
&&\quad +(35-15\sigma_1)\sigma_7+15\sigma_8)x^4+\dots
\end{eqnarray*}
with omitted terms divisible by $x^5$.

Section \ref{sectgeneralcase} discusses the combinatorics of
packings associated to an arbitrary (not necessarily generic) family
$\mathcal S_1,\dots,\mathcal S_n$ of subsets in a group. 

In Section \ref{sectexistence}, we specialize the 
results of Section \ref{sectgeneralcase} by applying them to 
generic packings. The underlying combinatorics are then simpler
and yield a proof of Proposition 
\ref{propexistence}, a crucial ingredient for establishing
the main result.

Sections \ref{sectpropfuncteq} and \ref{sectuniqness}
contain the proof of Propositions \ref{propfuncteqUformula} 
and \ref{propuniqueness} thus completing the proof of Theorem 
\ref{thmU}.

Section \ref{sectmoebius} uses Theorem \ref{thmU} and its proof 
for computing the M\"obius function of the poset of finite labelled 
hyperforests. An anonymous referee pointed out that this 
computation, a byproduct arising in the proof of our main result,
might be 
of independent interest. Remark \ref{thmweightedhypertrees} states
already known formulae for the enumeration of (weighted) labelled
hypertrees.

Section \ref{sectcomputexples} deals with computational aspects and examples.

Section \ref{sectasympt} contains a conjectural asymptotic formula
for the coefficients of $U(x,0,-1,-1,-1,\dots)$.

Section \ref{sectmodular} describes a few experimental 
observations concerning arithmetical properties of
the coefficients of $U(x,0,-1,-1,-1,\dots)$.

Section \ref{sectinttijn} is also experimental and describes a few integer
sequences related to the numbers $t_{i,j}(n)$ 
appearing, up to signs, as coefficients of the series $U$.

The paper ends with Section \ref{sectcov} discussing a few aspects 
of coverings which can be seen as dual objects of packings.


\section{Main result}\label{sectmain}

For $n=1,2,\dots$, we consider the following set $t_{i,j}(n)$ 
of strictly positive integers indexed by 
$i\in\{n+1,\dots,2n\}$ and $j\in\{0,1,\dots,2n-i\}$:
We set $t_{2,0}(1)=1$ and define $t_{i,j}(n)$ recursively by the
formula
\begin{eqnarray}\label{eqtriangle}
t_{i,j}(n)=(i-2)t_{i-1,j}(n-1)+t_{i-1,j-1}(n-1)+(i-3)t_{i-2,j}(n-1)
\end{eqnarray}
for $n\geq 2$. We set $t_{i,j}(n)=0$ in all other cases, 
i.e. if $i\leq n$ or $j< 0$ or $i+j>2n$.

Given a natural integer $n\geq 1$, 
the set of all ${n+1\choose 2}$ non-zero integers $t_{i,j}(n)$
can be organized into
a triangular array $T(n)$ with rows indexed by $\{n+1,\dots,2n\}$
and columns indexed by $\{0,\dots,n-1\}$ 
such that $T(n)$ determines $T(n+1)$ recursively
by Formula (\ref{eqtriangle}) reminiscent of the recurrence relation 
${n\choose k}={n-1\choose k-1}+{n-1\choose k}$ for binomial coefficients.
The first six triangular arrays $T(1),\dots,T(6)$ are
$$\begin{array}{r}
1\end{array}\qquad
\begin{array}{rr}
1&1\\1\end{array}\qquad
\begin{array}{rrr}
2&3&1\\5&3\\3\end{array}\qquad
\begin{array}{rrrrr}
6&11&6&1\\26&26&6\\
35&15\\
15\end{array}$$
$$\begin{array}{rrrrrrrr}
24&50&35&10&1\\
154&200&80&10\\
340&255&45\\
315&105\\
105\end{array}\qquad
\begin{array}{rrrrrrrrr}
120&274&225&85&15&1\\
1044&1604&855&190&15\\
3304&3325&1050&105\\
4900&2940&420\\
3465&945\\
945\end{array}$$
Observe that the first row of $T(1),T(2),\dots$ coincides, up to signs, 
with Stirling numbers of the first kind. More precisely, we have
\begin{eqnarray}
\label{firstrowstirling}
\sum_{k=0}^{n-1}t_{n+1,k}(n)x^{k+1}=\prod_{j=0}^{n-1}(x+j)=
(-1)^n\sum_{j=1}^nS_1(n,j)(-x)^j\  .
\end{eqnarray}
This is of course an easy 
consequence of the recurrence relation (\ref{eqtriangle}). 
The integers $t_{i,j}(n)$ seem to be related to
a few interesting integer-sequences, see Section 
\ref{sectinttijn} for examples.

We consider the formal power series $U\in A[[x]]$ with coefficients
in the ring $A=\mathbb Z[\sigma_1,\sigma_2,\sigma_3,\dots]$ of integral 
polynomials in $\sigma_1,\sigma_2,\dots$ defined by
\begin{eqnarray}\label{formulaU}
U(x,\sigma_1,\sigma_2,\dots)=
1-\sum_{n=1}^\infty x^n\sum_{i=n+1}^{2n}\sigma_i\sum_{j=0}^{2n-i}t_{i,j}(n)
(-\sigma_1)^j\ .
\end{eqnarray}

\begin{thm}\label{thmU} The number of packings of a generic family 
$\mathcal S_1,\dots,\mathcal S_n$ of $n$ non-empty subsets in a 
finite group $G$ with $N$ elements equals
\begin{eqnarray}\label{formulathmU}
N^nU(N^{-1},\sigma_1,\sigma_2,\dots)
\end{eqnarray}
for $U$ given by Formula (\ref{formulaU}) and for $\sigma_1,\sigma_2,\dots$
defined by
$$\sum_{j=0}^\infty \sigma_jt^j=\prod_{k=1}^n(1+\sharp(\mathcal S_k)t)\ .$$
\end{thm}

Remark that Formula (\ref{formulathmU}) of 
Theorem \ref{thmU} is polynomial of degree $n$ in $N$ 
for fixed complex numbers 
$\sigma_1,\sigma_2,\dots$ such that $\sigma_{n+1}=\sigma_{n+2}=\dots=0$.
Indeed, the coefficient of $x^m$ in $U(x,\sigma_1,\sigma_2,\dots)$
belongs to the ideal  generated by $\sigma_{m+1},\sigma_{m+2},
\dots,\sigma_{2m}$ of $\mathbb Z[\sigma_1,\sigma_2,\dots]$
and is thus zero for $m\geq n$ if $\sigma_{n+1}=\sigma_{n+2}=\dots
=0$.

The ingredients for proving Theorem \ref{thmU} are the following four
results:

\begin{prop} \label{propexistence} 
There exists a series $U\in\mathbb Z[[x,\sigma_1,\sigma_2,\dots]]$
such that Formula (\ref{formulathmU}) with $\sigma_1,\sigma_2,\dots$ 
defined as in Theorem \ref{thmU}
gives the number of packings
for every generic family of $n$ non-empty subsets in a finite group
with $N$ elements.

Moreover, the coefficient of a non-constant monomial 
$x^m$ in this series $U$ is of degree at most $2m$ with respect to the 
grading $\deg\sigma_i=i$ and belongs to the ideal of 
$\mathbb Z[\sigma_1,\sigma_2,\dots]$ generated by $\sigma_{m+1},
\sigma_{m+2},\dots,\sigma_{2m}$.
\end{prop} 

The proof of Proposition \ref{propexistence} relies on 
combinatorial properties
of intersection graphs encoding non-trivial intersections among subsets
$g_1\mathcal S_1,\dots,g_n\mathcal S_n$ of a group $G$. These properties
are encoded by the poset $\mathcal{HF}(n)$ of hyperforests with $n$
labelled vertices and order relation given by $F'\leq F$ if every 
hyperedge of $F'$ is contained in some hyperedge of $F$. The poset
$\mathcal{HF}(n)$ is a lattice with minimal element the trivial
graph defined by $n$ isolated labelled vertices and with maximal element
the hypertree consisting of a unique hyperedge containing all 
$n$ labelled vertices. Our proof of Proposition \ref{propexistence}
uses M\"obius inversion in 
$\mathcal HF(n)$. It needs only the existence (which is obvious) of a  
M\"obius function on the poset $\mathcal{HF}(n)$. The explicit 
description of $U$ given by Theorem \ref{thmU} allows however a
posteriori the computation (given by Proposition \ref{propmoebius})
of the M\"obius function of $\mathcal{HF}(n)$. Remark that the poset
$\mathcal{HT}_n$ of hypertrees with $n$ labelled vertices appearing for example 
in \cite{MM} is a subposet of the order dual of $\mathcal{HF}(n)$ obtained
by restricting the inverse order of $\mathcal{HF}(n)$ to 
the subset of all hypertrees in $\mathcal{HF}(n)$.

\begin{prop} \label{propfunctionalequationexistence}
A series $U$ as in Proposition \ref{propexistence} satisfies the functional 
equation 
\begin{eqnarray}\label{functionaleqU}
(1-\sigma_1 x)U(x,\sigma_1,\sigma_2,\sigma_3,\dots)=
U(x,\tilde\sigma_1,\tilde \sigma_2,\tilde\sigma_3,\dots)
\end{eqnarray}
where $\tilde \sigma_i=\sigma_{i-1}+\sigma_i$, using the convention
$\sigma_0=1$.
\end{prop}

{\bf Proof}
Equation (\ref{functionaleqU}) corresponds to equation 
(\ref{formeqalpha}) if $\sigma_1,\sigma_2,\dots$ are
elementary symmetric functions of a finite set of natural integers.
The general case follows by remarking that the algebra of 
symmetric polynomials is a free polynomial algebra on 
the set of elementary symmetric polynomials.\hfill$\Box$

\begin{prop} \label{propfuncteqUformula} The series $U$ defined by Formula
(\ref{formulaU}) satisfies the functional equation 
(\ref{functionaleqU}).
\end{prop}

\begin{prop} \label{propuniqueness} The functional equation (\ref{functionaleqU})
has at most one solution of the form $U=1+\dots$ such that the coefficient
of a nonconstant monomial $x^n$ is of degree at most $2n$ (with respect 
to the grading $\deg \sigma_i=i$) and belongs to the ideal generated by 
$\sigma_{n+1},\sigma_{n+2},\dots,\sigma_{2n}$ in $\mathbb Z[\sigma_1,\sigma_2,
\dots]$.
\end{prop}

{\bf Proof of Theorem \ref{thmU}} Proposition \ref{propexistence}
ensures the existence of a series enumerating packings 
of generic families in finite groups.
This series coincides with the series given by Formula 
(\ref{formulaU}) by Propositions \ref{propfunctionalequationexistence},
\ref{propfuncteqUformula} and \ref{propuniqueness}.
\hfill$\Box$

\begin{rem} Iterating identity (\ref{functionaleqU}) $n$ times 
we have
\begin{eqnarray*}
U(x,\sigma_1,\sigma_2,\dots)\prod_{j=0}^{n-1}(1-(\sigma_1+j)x)=
U(x,\tilde\sigma_1,\tilde\sigma_2,\tilde \sigma_3,\dots)
\end{eqnarray*}
where 
$$\tilde\sigma_k=\sum_{j=0}^{\min (k,n)}{n\choose j}\sigma_{k-j}\ .$$
A particular case is the specialization
$$U\left(x,{n\choose 1},{n\choose 2},{n\choose 3},\dots\right)=
\prod_{j=1}^{n-1}(1-jx)$$
associated to generic families $\mathcal S_1,\dots,\mathcal S_n$
given by $n$ singletons.
\end{rem}

\begin{rem} It is widespread lore that interesting combinatorial identities
have $q-$analogues generally encoding an additional feature
of the involved combinatorial objects. 
I do not know if the integers $t_{i,j}(n)$ or the series 
$U$ have such a $q-$analogue with interesting properties. 
\end{rem}


\section{Combinatorics of packings for arbitrary families 
$\mathcal S_1,\dots,\mathcal S_n$ of subsets in a group $G$}
\label{sectgeneralcase}

\subsection{Proof of Proposition \ref{propalphaineq}} 

A packing of $\mathcal S_1,\dots,\mathcal S_n$ given by $(g_1,\dots,g_n)
\in G^n$ extends to a packing 
$(g_1,\dots,g_n,g_{n+1})\in G^{n+1}$ of
$\mathcal S_1,\dots,\mathcal S_{n+1}$ if and only if 
$g_{n+1}\in G\setminus\left(\cup_{i=1}^n g_i\mathcal S_i
(\mathcal S_{n+1})^{-1}\right)$ where $\mathcal S^{-1}=\{g^{-1}\ \vert
\ g\in\mathcal S\}$.
Since $g_i\mathcal S_i(\mathcal S_{n+1})^{-1}$ contains at most $\sharp(
\mathcal S_{n+1})\sharp(\mathcal S_i)$ elements, we have the first 
inequality.

Considering a fixed element $h\in\mathcal S_{n+1}$ we have
the inequality
\begin{eqnarray*}
\sharp\left(\cup_{i=1}^ng_i\mathcal S_i(\mathcal S_{n+1})^{-1}\right)&\geq&
\sharp\left(\cup_{i=1}^ng_i\mathcal S_ih^{-1}\right)
=\sharp\left(\cup_{i=1}^ng_i\mathcal S_i\right)\ .
\end{eqnarray*}
For a packing $(g_1,\dots,g_n)$, we have 
$$\sharp\left(\cup_{i=1}^ng_i\mathcal S_i\right)=
\sum_{i=1}^n\sharp(\mathcal S_i)$$
showing the second inequality. 

Both inequalities are sharp if $\sharp(\mathcal S_{n+1})=1$. This proves
equality (\ref{formeqalpha}).
\hfill$\Box$

\subsection{Intersection graphs}\label{sectintersgraph}

We fix a group $G$ and a family $\mathcal S_1,\dots,\mathcal S_n$ 
of $n$ non-empty subsets in $G$. 
Given an element $\mathbf g=(g_1,\dots,g_n)$ of $G^n$,
we consider the corresponding
{\it intersection graph} $\mathcal I(\mathbf g)$ with vertices 
$1,\dots,n$ and edges 
$\{i,j\}$ between distinct vertices $i,j$ if 
$g_i\mathcal S_i\cap g_j\mathcal S_j\not=\emptyset $ in $G$. 
Observe that $\mathbf g=(g_1,\dots,g_n)$
in $G^n$ defines a packing if and only if 
$\mathcal I(\mathbf g)$ is the trivial graph with $n$ isolated vertices.

Given a finite simple graph $\Gamma$ with vertices $1,\dots,n$ 
and edges $E(\Gamma)$, we consider the set
$$\mathcal R_\Gamma=\{(g_1,\dots,g_n)\in G^n\ \vert\ 
g_i\mathcal S_i\cap g_j\mathcal S_j\not=\emptyset\hbox{ for every }
\{i,j\}\in E(\Gamma)\}\ .$$

An element $\mathbf g$ in $G^n$ belongs thus to $\mathcal R_\Gamma$ 
if and only if $\Gamma$ is 
a subgraph of the intersection graph $\mathcal I(\mathbf g)$.

We denote by $\mathcal E_\Gamma$
the set of equivalence classes of $\mathcal R_\Gamma$ defined by 
$(g_1,\dots,g_n)\sim(h_1,\dots,h_n)$ if 
$g_ih_i^{-1}=g_jh_j^{-1}$ for every edge $\{i,j\}$ of $\Gamma$.
Two elements $\mathbf g=(g_1,\dots,g_n)$ and $\mathbf h=(h_1,\dots,h_n)$ of
$\mathcal R_\Gamma$ represent thus the same
equivalence class of $\mathcal E_\Gamma$ 
if and only if the map $i\longmapsto g_ih_i^{-1}$
is constant on (vertices of) connected components.

\begin{prop} \label{propRGammaN} 
Suppose that $G$ is a finite group with $N$ elements.
We have then
$$\sharp(\mathcal R_\Gamma)=\sharp(\mathcal E_\Gamma)N^{c(\Gamma)}$$
where $c(\Gamma)$ denotes the number of connected components of $\Gamma$.
\end{prop} 

{\bf Proof} We set $c=c(\Gamma)$ and 
we denote the connected components of $\Gamma$ 
by $\Gamma_1,\dots,\Gamma_c$. We get a free action 
of $G^c$ on $\mathcal R_\Gamma$ by considering
$$(a_1,\dots,a_c)\cdot(g_1,\dots,g_n)
\longmapsto (a_{\gamma(1)}^{-1}g_1,\dots,
a_{\gamma(n)}^{-1}g_n)$$ 
where $\gamma(i)\in\{1,\dots,c\}$  is defined by 
the inclusion of the vertex $i$ in the $\gamma(i)-$th connected component
$\Gamma_{\gamma(i)}$ of $\Gamma$. Orbits in $\mathcal R_\Gamma$
of this action are thus in one-to-one correspondence with 
equivalence classes of $\mathcal E_\Gamma$.
\hfill$\Box$

\begin{rem}
The set $\mathcal E_\Gamma$ associated to a graph 
$\Gamma$ with $c$ connected components contains at most
$\left(\max_{i}\sharp(\mathcal S_i)\right)^{2n-2c}$ distinct 
equivalence classes. Indeed, we have 
$\mathcal R(\Gamma')\subset \mathcal R(\Gamma)$ if $\Gamma$ is a subgraph 
of $\Gamma'$. Replacing $\Gamma$ by a spanning forest, we can thus assume
that $\Gamma$ is a forest. The equivalence class of an element 
${\mathbf g}\in \mathcal R(\Gamma)$ is now determined by the 
relative positions of $g_i\mathcal S_i$ and $g_j\mathcal S_j$ for all 
$n-c$ edges $\{i,j\}$ of the forest $\Gamma$ and the number
of different relative positions of $g_i\mathcal S_i$ and $g_j\mathcal S_j$
is at most $\sharp(\mathcal S_i)\sharp(\mathcal S_j)\leq  
\left(\max_{i}\sharp(\mathcal S_i)\right)^2$.
\end{rem}

\subsection{M\"obius inversion}

\begin{prop}\label{propalphaN}
The number $\alpha=\alpha(G;\mathcal S_1,\dots,\mathcal S_n)$ of packings 
of a family $\mathcal S_1,\dots,\mathcal S_n$ in a finite group $G$ with $N$
elements is given by
$$\alpha=\sum_{\Gamma\in\mathcal B} 
(-1)^{e(\Gamma)}\sharp(\mathcal E_\Gamma)N^{c(\Gamma)}$$
where the sum is over the Boolean poset $\mathcal B$ of all 
$2^{n\choose 2}$ simple graphs with vertices $1,\dots,n$ and where 
$e(\Gamma)=\sharp(E(\Gamma))$, 
respectively $c(\Gamma)$, denotes the number of edges,
respectively connected components, of a graph $\Gamma\in\mathcal B$.
\end{prop}

{\bf Proof} Proposition \ref{propRGammaN} shows that it is enough 
to prove the equality
$$\alpha=\sum_{\Gamma\in\mathcal B}(-1)^{e(\Gamma)}
\sharp(\mathcal R_\Gamma)\ .$$
An element $\mathbf g=(g_1,\dots,g_n)\in G^n$
defines a packing if and only if its intersection graph 
$\mathcal I(\mathbf g)$ is trivial. It provides thus a contribution of $1$ to
$\alpha$ in this case since it is only involved as an element 
of $\mathcal R_\Gamma$ if
$\Gamma$ is the trivial graph with isolated vertices $1,\dots,n$ and 
no edges.

An element $\mathbf g=(g_1,\dots,g_n)\in G^n$
with non-trivial intersection graph $\mathcal I(\mathbf g)$
containing $e\geq 1$ edges yields a contribution of $0$ to $\alpha$ 
since contributions coming from 
the $2^{e-1}$ subgraphs of $\mathcal I(\mathbf g)$ containing an even 
number of edges cancel out with contributions associated to the 
$2^{e-1}$ subgraphs 
having an odd number of edges.
\hfill$\Box$

\begin{rem} Introducing 
$$\alpha_\Gamma=\{{\mathbf g}\in G^n\ 
\vert\ \mathcal I({\mathbf g})=\Gamma\}\ ,$$
we have $\alpha=\alpha_T$ where $T$ denotes the trivial graph 
with $n$ isolated vertices $1,\dots,n$.
Our proof of 
Proposition \ref{propalphaN} computes $\alpha$ by applying M\"obius inversion
(more precisely, its dual form, see Proposition 3.7.2 of \cite{St}) 
\begin{eqnarray}\label{formulamoebiusgeneral}
\alpha=\sum_{\Gamma\in\mathcal B} \mu(\Gamma)\sharp(
\mathcal R_\Gamma)
\end{eqnarray}
(with $\mu(\Gamma)=(-1)^{e(\Gamma)}$ denoting the M\"obius function
of the Boolean lattice $\mathcal B$ of all simple graphs on $1,\dots,n$)
to the numbers 
$$\sharp(\mathcal R_\Gamma)=
\sum_{\Gamma'\supset\Gamma}\alpha_{\Gamma'}$$
given by Proposition \ref{propRGammaN}.
\end{rem}


\section{Proof of Proposition \ref{propexistence}: 
Combinatorics of generic packings}\label{sectexistence}

A \emph{hypergraph} consists of a set $\mathcal V$ of vertices
and of a set of hyperedges where a hyperedge is a subset of 
$\mathcal V$ containing at least $2$ vertices. Two vertices 
are \emph{adjacent} if they belong to
a common hyperedge. A \emph{path} is a
sequence of consecutively adjacent vertices.
A hypergraph is \emph{connected} if any pair of vertices can
be joined by a path.
A \emph{cycle} is a closed path involving only distinct vertices.
A \emph{hyperforest} is a hypergraph with
distinct hyperedges intersecting in at most a common vertex and
with every cycle contained in a hyperedge. A
\emph{hypertree} is a connected hyperforest.

The \emph{primal graph} of a hypergraph with vertices 
$\mathcal V$ is the ordinary graph with vertices $\mathcal V$ 
and ordinary edges encoding adjacency in the hypergraph.
An ordinary graph $\Gamma$ is the primal graph of a hyperforest if and only
if every cycle and every
edge of $\Gamma$ is contained in a unique maximal complete
subgraph. Maximal complete subgraphs of such a graph 
$\Gamma$ are in one-to-one
correspondence with hyperedges of the associated hyperforest.
Primal graphs of hyperforests are often called \emph{block-graphs}
or \emph{chordal and diamond-free graphs}.
In the sequel, we identify generally hyperforests with their primal 
graphs.

\begin{lem} \label{lemedge}
The intersection 
$g_i\mathcal S_i\cap g_j\mathcal S_j$ associated 
to an edge $\{i,j\}$ in an intersection graph 
$\mathcal I(\mathbf g)$ is reduced to a unique element
if $\mathcal S_1,\dots,\mathcal S_n$ is a generic family of $G$.
\end{lem}

{\bf Proof} Otherwise there exist two distinct elements $a_i,b_i\in
\mathcal S_j$ and two distinct elements $a_j,b_j\in\mathcal S_j$ such that
$g_ia_i=g_jb_j$ and $g_ja_j=g_ib_i$. This shows 
$$g_ia_ib_j^{-1}g_j^{-1}g_ja_jb_i^{-1}g_i^{-1}=e$$
and implies the relation
$b_i^{-1}a_ib_j^{-1}a_j=e$ with $b_i^{-1}a_i\in\mathcal S_i^{-1}\mathcal S_i
\setminus\{e\}$ and $b_j^{-1}a_j\in\mathcal S_j^{-1}\mathcal S_j
\setminus\{e\}$ in contradiction with genericity of the family
$\mathcal S_1,\dots,\mathcal S_n$. \hfill$\Box$

\begin{prop} \label{propinterblock}
Intersection graphs of generic families are (primal graphs of) hyperforests.
\end{prop}

{\bf Proof} Consider $k$ cyclically consecutive vertices 
$i_1,i_2,\dots,i_{k-1},i_k,i_{k+1}=i_1$ in an intersection graph
$\mathcal I(\mathbf g)$ of a generic family 
$\mathcal S_1,\dots,\mathcal S_n\subset G$.
Lemma \ref{lemedge} implies the existence of unique elements 
$a_{i_j}\in \mathcal S_{i_j}$ and $b_{i_{j+1}}\in\mathcal S_{i_{j+1}}$
such that $g_{i_j}a_{i_j}=g_{i_{j+1}}b_{i_{j+1}}$ for every edge 
$\{i_j,i_{j+1}\}$ of $C$. We get thus the relation
\begin{eqnarray*}
g_{i_1}a_{i_1}(g_{i_2}b_{i_2})^{-1}g_{i_2}a_{i_2}(g_{i_3}b_{i_3})^{-1}\cdots
g_{i_k}a_{i_k}(g_{i_1}b_{i_1})^{-1}=e
\end{eqnarray*}
which is conjugate to the relation
\begin{eqnarray*}
\left(b_{i_1}^{-1}a_{i_1}\right)\left(b_{i_2}^{-1}a_{i_2}\right)\cdots
\left(b_{i_k}^{-1}a_{i_k}\right)=e\ .
\end{eqnarray*}
Genericity of the family $\mathcal S_1,\dots,\mathcal S_n$
implies $a_{i_j}=b_{i_j}$ for all $j$. The sets
$g_{i_j}\mathcal S_{i_j}$ intersect thus in the common element
$g_{i_1}a_{i_1}=\dots=g_{i_k}a_{i_k}$ (which is the unique
common element of pairwise distinct sets in 
$\{g_{i_1}\mathcal S_{i_1},\dots,g_{i_k}\mathcal S_{i_k}\}$
by Lemma \ref{lemedge}).
All elements $i_1,\dots,i_k$ of $\mathcal I(\mathbf g)$
are thus adjacent vertices contained in a common maximal complete
subgraph of $\mathcal I(\mathbf g)$. 

Suppose now that an edge $\{i,j\}$ belongs to two distinct maximal 
complete subgraphs $K$ and $K'$ of $\mathcal I({\mathbf g})$.
Maximality of $K$ and $K'$ implies the existence of vertices
$k\in K\setminus K'$ and $k'\in K'\setminus K$. 
Thus we get triplets of mutually adjacent vertices $i,j,k\subset K$ and
$i,j,k'\subset K'$. Lemma \ref{lemedge} shows that $g_i\mathcal S_i
\cap g_j\mathcal S_j$ is reduced to a unique element $a$. We have thus 
$g_i\mathcal S_i\cap g_j\mathcal S_j\cap g_k\mathcal S_k=\{a\}
\subset K$. Similarly, we get $a\in g_{k'}\mathcal S_{k'}$.
This implies $k'\in K$ in contradiction with $k'\in K'\setminus K$.

Distinct maximal complete subgraphs of $\mathcal I({\mathbf g})$ intersect 
thus at most in a common vertex and 
every cycle of $\mathcal I({\mathbf g})$ is contained in a unique
maximal complete subgraph of $\mathcal I({\mathbf g})$.
This implies that $\mathcal I(\mathbf g)$
is (the primal graph of) a hyperforest.\hfill$\Box$

For the sake of concision, we identify in the sequel such an intersection
graph $\mathcal I(\mathbf g)$ with the corresponding hyperforest.

Applying the proof of Proposition \ref{propinterblock}
to a Hamiltonian cycle  visiting all vertices of 
a hyperedge $\{i_1,\dots,i_k\}$ in 
an intersection graph $\mathcal I(\mathbf g)$ associated to a
generic family we get the following result: 

\begin{prop} Given a hyperedge $\{i_1,\dots,i_k\}$ in the intersection graph 
$\mathcal I(\mathbf g)$ of a generic family $\mathcal S_1,\dots,\mathcal S_n$,
there exists a unique element $a\in G$ such that 
$g_{i_l}\mathcal S_{i_l}\cap g_{i_m}\mathcal S_{i_m}=\{a\}$
for every pair of distinct vertices $i_l,i_m$ in $\{i_1,\dots,i_k\}$.
\end{prop}

\begin{prop}\label{propcardCCeq} 
Let $\Gamma$ be a hyperforest with vertices $1,\dots,n$
indexing the subsets $\mathcal S_1,\dots,\mathcal S_n$
of a generic family in a group $G$. Defining the equivalence
relation $\mathcal E_F$ of a hyperforest $F$ as in Section 
\ref{sectintersgraph}, we have 
$$\sharp(\mathcal E_F)=\prod_{j=1}^n\left(\sharp(\mathcal S_j)\right)^{
\deg_F(j)}$$
with $\deg_F(j)$ denoting the degree of $j$ defined as the 
number of distinct hyperedges containing the vertex $j$.
\end{prop}

{\bf Proof} Let $e=\{i_1,\dots,i_k\}$ be a hyperedge of an 
intersection graph $\mathcal I({\mathbf g})$. Since 
$\bigcap_{j=1}^kg_{i_j}\mathcal S_{i_j}$ is reduced to a unique 
element $a_e\in G$, we get a map $\mu_e:\{i_1,\dots,i_k\}\longrightarrow
G$ such that $\mu_e(i_j)\in\mathcal S_{i_j}$
by setting $\mu_e(i_j)=g_{i_j}^{-1}a_e$. This map depends only on the
equivalence class in $\mathcal E_{\mathcal I({\mathbf g})}$ 
of $\mathcal I({\mathbf g})$ and the set of all such maps
determines the equivalence class of $\mathcal I({\mathbf g})$
in $\mathcal E_F$ for any hyperforest $F$ contained in 
$\mathcal I({\mathbf g})$. Since all cycles of a hyperforest
are contained in hyperedges, all possible choices of the 
maps $\mu_e$ associated 
to hyperedges of $F$ correspond to equivalence classes of $\mathcal E_F$. 
Different choices yield inequivalent classes. 
The set $\mathcal E_F$ of all equivalence classes is thus in 
one-to-one correspondence with the set $\prod_{j=1}^n
\mathcal S_j^{\mathop{\deg_F(j)}}$.
\hfill$\Box$

{\bf Proof of Proposition \ref{propexistence}} Setting $s_i=
\sharp(\mathcal S_i)$, Proposition \ref{propcardCCeq}
can be rewritten as the identity 
$$\sharp(\mathcal E_F)=\prod_{j=1}s_j^{\deg_F(j)}$$
for every hyperforest $F$ with vertices $\{1,\dots,n\}$.
We denote by $\mathcal{HF}(n)$ the set of all hyperforests with vertices 
$\{1,\dots,n\}$. The set $\mathcal{HF}(n)$ is partially ordered by 
inclusion by setting $F'\leq F$ for $F',F\in\mathcal{HF}(n)$ if
every hyperedge of $F'$ is contained in some hyperedge of $F$.
Equivalently, $F'\leq F$ if adjacent vertices
of $F'$ are also always adjacent in $F$. The primal graph underlying $F'$ is
thus a subgraph of the primal graph underlying $F$ if $F'\leq F$.
Denoting by $\mu$ the M\"obius function of the poset $\mathcal{HF}(n)$, 
the number $\alpha=\alpha(G;\mathcal S_1,\dots,
\mathcal S_n)$ of 
packings of a generic family $\mathcal S_1,\dots,\mathcal S_n$
in a group $G$ of order $N$ is given by 
\begin{eqnarray}\label{formaproof}
\alpha=\sum_{F\in\mathcal{HF}(n)}\mu(F)N^{c(F)}\prod_{j=1}^n
s_j^{\deg_F(j)}
\end{eqnarray}
(with $c(F)$ denoting the number of connected components of 
a hyperforest $F$). Since the M\"obius function of $\mathcal{HF}(n+1)$
restricts to the M\"obius function of $\mathcal{HF}(n)$, the summation
over $\mathcal{HF}(n)$ in Formula (\ref{formaproof}) can be extended
(after setting $s_i=0$ for $i>n$ and using the convention $0^0=1$)
over the poset $\mathcal{HF}$ of all hyperforests with vertices 
$\mathbb N\setminus\{0\}$ such that almost all vertices
are isolated (only
finitely many vertices have strictly positive degree).

Summing over all possible labellings of an unlabelled hyperforest
and remarking that the M\"obius function is invariant under permutations
of labels
shows that $\alpha$ is a symmetric function of $s_1,\dots,s_n$.
This expresses $\alpha$ as a polynomial in 
$\sigma_1=\sum s_i,\sigma_2=\sum_{i<j} s_is_j,\dots$ 
with contributions coming from finite unlabelled hyperforests.

More precisely, contributions to the
coefficient $N^{n-m}$ of $\alpha$ given by 
Formula (\ref{formaproof}) come from 
hyperforests with vertices $\{1,\dots,n\}$
consisting of $c\geq 0$ non-trivial hypertrees
involving $m+c\leq n$ vertices of strictly positive degrees
and $n-m-c$ isolated vertices. The equalities 
$s_{n+1}=s_{n+2}=\dots=0$ imply that the summation over 
$\mathcal{HF}(n)$ in (\ref{formaproof}) can be extended to a summation
over all hyperforests in $\mathcal{HF}(l)$ for an arbitrary integer
$l\geq n$ since the vertices $n+1,n+2,\dots$ have to be isolated vertices
of a hyperforest yielding a non-zero contribution to $\alpha$. 
Observe now that we have $c\leq m$ since every non-trivial hypertree 
contains at least two vertices. For a fixed value of $m$, a
labelled hyperforest with non-zero contribution to $\alpha$
has thus at most $2m$ non-isolated vertices (with equality achieved by 
a hyperforest consisting of $m$ isolated edges
joining $2m$ distinct vertices).
Summing over unlabelled hyperforests and considering the associated 
symmetric functions $\sigma_1,\sigma_2,\dots$ in $s_1,s_2,\dots$
(obtained by a summation over all possible distinct labellings of the
underlying unlabelled hypertrees) we see that 
the contribution associated to an 
unlabelled hyperforest with $m+c\leq 2m$ non-isolated vertices
is in the ideal of $\mathbb Z[x,\sigma_1,\sigma_2,\dots]$
generated by $\sigma_{m+c},\sigma_{m+c+1},\sigma_{m+c+2},\dots$.
The degree in $\sigma_1,\sigma_2,\dots$ (with respect to the 
grading $\deg(\sigma_i)=i$) of such a contribution
is maximal and equals $2m$ for ordinary unlabelled forests having $m$
ordinary edges. Indeed, let $F$ be a hyperforest with $c$ connected
components and $m+c$ vertices of strictly positive degrees. 
Replacing a hyperedge $E$ of $F$ involving $k\geq 3$ vertices by
a tree consisting of
$k-1$ ordinary edges connecting all vertices of $E$
increases the degree-sum of all vertices by $k-2>0$ and yields
a contribution of higher degree. Contributions of maximal 
degree correspond thus to ordinary forests on 
$m+c$ vertices with $c$ connected components. Such a forest has 
$m$ edges and yields a contribution of degree $2m$ 
with respect to the grading $\deg(\sigma_i)=i$.
This ends the proof of Proposition \ref{propexistence}.
\hfill$\Box$


\section{Proof of Proposition 
\ref{propfuncteqUformula}}\label{sectpropfuncteq}

We have to show that 
$$U=
1-\sum_{n=1}^\infty x^n\sum_{i=n+1}^{2n}\sigma_i\sum_{j=0}^{2n-i}t_{i,j}(n)
(-\sigma_1)^j$$
defined by Formula (\ref{formulaU}) satisfies the functional equation 
$$(1-\sigma_1 x)U(x,\sigma_1,\sigma_2,\sigma_3,\dots)=
U(x,1+\sigma_1,\sigma_1+\sigma_2,\sigma_2+\sigma_3,\dots)\ ,$$ 
see (\ref{functionaleqU}). 
Both sides of (\ref{functionaleqU}) have the same constant term $1$ and
involve only non-constant monomials of the form 
$\sigma_i\sigma_1^jx^n$. It is thus enough to check 
that coefficients of both sides of (\ref{functionaleqU}) 
agree for such monomials. This is easily checked for the coefficient of $x$. 
For a general monomial of the form $\sigma_i\sigma_1^jx^n$,
equation (\ref{functionaleqU}) amounts to the identity 
\begin{eqnarray*}
&&-\left((-1)^jt_{i,j}(n)-(-1)^{j-1}t_{i,j-1}(n-1)\right)\\
&=&-\left(\sum_{k=j}^{2n-i}t_{i,k}(n)(-1)^k{k\choose j}
+\sum_{k=j}^{2n-i-1}t_{i+1,k}(n)(-1)^k{k\choose j}\right)
\end{eqnarray*}
or equivalently to
\begin{eqnarray}\label{idrecurs}
t_{i,j}(n)+t_{i,j-1}(n-1)=\sum_k(-1)^{k+j}{k\choose j}\left(
t_{i,k}(n)+t_{i+1,k}(n)\right)
\end{eqnarray}
where $\sum_kf(k)=\sum_{k\in\mathbb Z}f(k)$ since ${k\choose j}\left(
t_{i,k}(n)+t_{i+1,k}(n)\right)=0$ for $k<j$ or $k> 2n-i$. 
We prove (\ref{idrecurs}) by induction on $n$. A straightforward computation
shows that it holds $n=2$.
Applying the recursion relation (\ref{eqtriangle}) which holds for all
$i,j\in\mathbb Z$ if $n\geq 2$ to 
the right side 
$$R=\sum_k(-1)^{k+j}{k\choose j}(t_{i,k}(n)+t_{i+1,k}(n))$$
of (\ref{idrecurs})
we get
\begin{eqnarray*}
R
&=&\sum_k(-1)^{k+j}{k\choose j}\Big(\\
&&(i-2)t_{i-1,k}(n-1)+t_{i-1,k-1}(n-1)+(i-3)t_{i-2,k}(n-1)\\
&&\quad +(i-1)t_{i,k}(n-1)+t_{i,k-1}(n-1)+(i-2)t_{i-1,k}(n-1)\Big)\\
&=&L+C\end{eqnarray*}
where
\begin{eqnarray*}
L&=&(i-2)\sum_k(-1)^{k+j}{k\choose j}(t_{i-1,k}(n-1)+t_{i,k}(n-1))\\
&&\quad +\sum_{k} (-1)^{k+j-1}{k\choose j-1}(t_{i-1,k}(n-1)+t_{i,k}(n-1))\\
&&\quad +(i-3)\sum_k(-1)^{k+j}{k\choose j}(t_{i-2,k}(n-1)+t_{i-1,k}(n-1))
\end{eqnarray*}
and
\begin{eqnarray*}
C&=&-\sum_k(-1)^{k+j-1}{k\choose j-1}(t_{i-1,k}(n-1)+t_{i,k}(n-1))\\
&&\quad +\sum_{k}(-1)^{k+j}{k\choose j}(t_{i-1,k-1}(n-1)+t_{i,k-1}(n-1))\\
&&\quad +\sum_{k}(-1)^{k+j}{k\choose j}(t_{i,k}(n-1)+t_{i-1,k}(n-1))\\
&=&\quad \sum_k(-1)^{k+j}{k\choose j-1}(t_{i-1,k}(n-1)+t_{i,k}(n-1))\\
&&\quad -\sum_{k}(-1)^{k+j}{k+1\choose j}(t_{i-1,k}(n-1)+t_{i,k}(n-1))\\
&&\quad +\sum_{k}(-1)^{k+j}{k\choose j}(t_{i-1,k}(n-1)+t_{i,k}(n-1))\\
&=&\sum_k(-1)^{k+j}\left({k\choose j-1}-{k+1\choose j}+{k\choose j}\right)
 (t_{i-1,k}(n-1)+t_{i,k}(n-1))
\end{eqnarray*}
which shows $C=0$ since ${k+1\choose j}={k\choose j-1}+{k\choose j}$.

Using induction on $n$ and applying (\ref{idrecurs}) we get
\begin{eqnarray*}
L&=&(i-2)(t_{i-1,j}(n-1)+t_{i-1,j-1}(n-2))\\
&&\quad +(t_{i-1,j-1}(n-1)+t_{i-1,j-2}(n-2))\\
&&\quad +(i-3)(t_{i-2,j}(n-1)+t_{i-2,j-1}(n-2))
\end{eqnarray*}
We have thus 
\begin{eqnarray*}
L&=&(i-2)t_{i-1,j}(n-1)+t_{i-1,j-1}(n-1)+(i-3)t_{i-2,j}(n-1)\\
&&+(i-2)t_{i-1,j-1}(n-2)+t_{i-1,j-2}(n-2)+(i-3)t_{i-2,j-1}(n-2)
\end{eqnarray*}
and applying (\ref{eqtriangle}) we get 
$$L=t_{i,j}(n)+t_{i,j-1}(n-1)$$
which is the left side of (\ref{idrecurs}). \hfill$\Box$


\section{Proof of Proposition \ref{propuniqueness}}\label{sectuniqness}

Assuming the existence of two distinct series $U_1,U_2$ 
fulfilling the requirements of 
Proposition \ref{propuniqueness}, the difference 
$D=U_1-U_2=\sum_{n=1}^\infty D_nx^n$ satisfies all hypotheses 
except for the value of its constant term.
Since $U_1$ and $U_2$ are different, there exists a minimal natural 
integer $n\geq 1$ such that $D_n\not=0$.
Let $m\geq n+1$ be the smallest integer such that 
$D_n=\sum_{k=m}^{2n}\sigma_kC_{n,k}$ with
$C_{n,k}\in\mathbb C[\sigma_1,\sigma_2,\dots]$ and 
$C_{n,m}\not=0$. Since $D_n$ is of degree $\leq 2n$ with respect to
the grading given by $\deg(\sigma_i)=i$, we have 
$C_{n,m}\in \mathbb C[\sigma_1,\dots,\sigma_{2n-m}]\subset 
\mathbb C[\sigma_1,\dots,\sigma_{n-1}]$.

Equation (\ref{functionaleqU}) and minimality of $n$ imply
$$D_n(1+\sigma_1,\sigma_1+\sigma_2,\sigma_2+\sigma_3,\dots)=
D_n(\sigma_1,\sigma_2,\sigma_3,\dots)$$
or equivalently
$$\sum_{k=m}^{2n}(\sigma_{k-1}+\sigma_k)C_{n,k}(1+\sigma_1,
\sigma_1+\sigma_2,\sigma_2+\sigma_3,\dots)=
\sum_{k=m}^{2n}\sigma_kC_{n,k}(\sigma_1,\sigma_2,\sigma_3,\dots)\ .$$
Comparison of both sides modulo the ideal $I$ generated by 
$\sigma_m,\sigma_{m+1},\sigma_{m+2},\dots$
gives 
$$C_{n,m}(1+\sigma_1,\sigma_1+\sigma_2,\sigma_2+\sigma_3,
\dots)=0\ .$$
Algebraic independency of the symmetric functions $\sigma_1,\sigma_2,\dots$
shows thus $C_{n,m}=0$ in contradiction with our assumption.\hfill$\Box$.


\section{The M\"obius function for the poset of finite 
labelled hyperforests}\label{sectmoebius}

Let $\mathcal P$ be a poset (partially ordered set) such that 
$\mathcal P$ has a unique minimal element $m$ and
$\{y\in\mathcal P\ \vert\ y<x\}$ is finite for all $x\in\mathcal P$.
This allows the recursive definition of a M\"obius function $\mu$ by
setting $\mu(m)=1$ and $\mu(x)=-\sum_{y<x}\mu(y)$ for all $x>m$.
Given a function $f:\mathcal P\longrightarrow \mathbb C$ 
with finite support, the value $f(m)$
can then be recovered from the function 
$g(x)=\sum_{y\geq x}f(y)$ using M\"obius inversion
$$f(m)=\sum_{x\in\mathcal P} \mu(x)g(x)\ ,$$
see Proposition 3.7.2 of \cite{St} (we use only the values $\mu(m,x)$
of the M\"obius function and write $\mu(x)=\mu(m,x)$ in analogy with the
usual, well-known number-theoretic M\"obius function of natural integers). 
M\"obius inversion was the main ingredient in the proof of 
Proposition \ref{propexistence}.
The poset $\mathcal{HF}$ of hyperforests
consisting of all hyperforests (ordered by inclusion) 
with finitely many hyperedges and vertices 
$1,2,3,4,\dots$ has a minimal element given by the trivial 
graph having only isolated vertices. The number
$$\sharp\lbrace F'\in \mathcal{HF}\ \vert\ F'\subset F\}$$
of all hyperforests contained in a given hyperforest 
$F\in\mathcal{HF}$ with $n$ hyperedges of degrees
$d_1,\dots,d_n$ is bounded by the number 
$2^{{d_1\choose 2}+\dots+{d_n\choose 2}}$ of (ordinary) subgraphs of
the primal graph underlying $F$.
The poset $\mathcal{HF}$ has thus a M\"obius function.

\begin{prop} \label{propmoebius}
The M\"obius function $\mu(F)$ of a hyperforest
$F$ in the poset $\mathcal{HF}$ of
all vertex-labelled hyperforests with finitely many hyperedges is given by
\begin{eqnarray}\label{formmoebforest}
\mu(F)=\prod_{j\geq 2}(-(j-2)!)^{\kappa_j}
\end{eqnarray}
where $\kappa_j$ denotes the number of hyperedges
involving exactly $j$ vertices of $F$.
\end{prop}

\begin{rem} \label{remposethyper}
The poset $\mathcal{HF}$ is in fact a lattice with wedge 
$F_1\wedge F_2$ given by the intersection and join $F_1
\vee F_2$ given by the smallest hyperforest containing $F_1$
and $F_2$.
\end{rem}

{\bf Proof of Proposition \ref{propmoebius}} 
Remark first that the order relation induced on 
subforests of a given hyperforest $F\in\mathcal{HF}$ 
is the product order of all 
order-relations on hyperedges of $F$.
An easy argument (or Proposition 3.8.2 of \cite{St}) shows thus that we have
$$\mu(F)=\prod_{e\in E(F)}\mu(e)$$
where $E(F)$ denotes the set of hyperedges of 
$F$ and where $\mu(e)$ 
is the M\"obius function  restricted to a hyperedge $e\in E(F)$.
This can of course be rewritten as
$$\mu(F)=\prod_{j\geq 2}\mu(K_j)^{\kappa_j}$$
where $K_j$ is an abitrary hyperedge on $j$ 
labelled vertices and where 
$\kappa_j$ is the number of hyperedges having $j$ vertices of $F$.

The proof of Proposition \ref{propexistence} shows that 
$\mu(K_j)$ coincides with the coefficient of $\sigma_{j+1}x^j$
in $U$. By Theorem \ref{thmU} (whose proof needs only the existence
but not the exact determination of the M\"obius function),
this coefficient equals 
$-t_{j+1,0}(j)=-(j-2)!$ where the last identity follows easily 
from Formula (\ref{eqtriangle}) defining the integers
$t_{i,j}(n)$  recursively.
\hfill$\Box$

\begin{rem} It would be interesting to have a simple direct proof 
that $\mu(K_n)=-(n-2)!$ for a hypergraph $K_n\in \mathcal{HF}$ 
reduced to a unique hyperedge involving $n\geq 2$ vertices.
\end{rem}

\begin{rem}\label{thmweightedhypertrees}
Let $\mathcal{HT}_k(n)$ be the finite set of all hypertrees with 
$k$ hyperedges and $n$ labelled vertices. Denoting by $\sharp(e)$ the 
number of vertices involved in a hyperedge $e$, we have
\begin{eqnarray}\label{preciseweightedformula}
\sum_{T\in\mathcal{HT}_k(n)}\prod_{e\in\mathcal E(T)}
(\sharp(e)-2)!\prod_{j=1}^n s_j^{\deg(j)}=(-1)^{n+k+1}
\sigma_n\sigma_1^{k-1}S_1(n-1,k)
\end{eqnarray}
where $\sigma_1=\sum_{j=1}^n s_j$ and $\sigma_n=\prod_{j=1}^n s_j$
and where $S_1(n,k)$
denotes the Stirling number of the first kind defined by 
$$\sum_{k=0}^nS_1(n,k)x^k=x(x-1)(x-2)\cdots (x-n+1)=
\prod_{j=0}^{n-1}(x-j)\ .$$
Indeed, the proof of Proposition \ref{propexistence}
shows that a hyperforest with $n$ non-isolated vertices,
$k$ hyperedges and $c$ connected components yields only 
contributions to the coefficients of $x^{n-c}\sigma_{n+s}\sigma_1^{k-1-s}$
for $s=0,\dots,k-1$. The coefficient 
of $x^{n-1}\sigma_n\sigma_1^{k-1}$ in $U$ is thus obtained from
contributions from all elements in the set $\mathcal{HT}_k(n)$
of hypertrees with $n$ non-isolated vertices $\{1,\dots,n\}$ and 
$k$ hyperedges. This coefficient equals $(-1)^{n+1}\sigma_n\sigma_1^{k-1}
S_1(n-1,k)$ by Formulae (\ref{formulaU}) and (\ref{firstrowstirling}).
Formulae (\ref{formaproof}) and (\ref{formmoebforest}) show that 
a hypertree $T\in \mathcal{HT}_k(n)$ contributes a summand
given by $(-1)^k\prod_{e\in\mathcal E(T)}
(\sharp(e)-2)!\prod_{j=1}^n s_j^{\deg(j)}$ to the coefficient
of $x^{n-1}\sigma_n\sigma_1^{k-1}$ in $U$.

Setting $s_1=\dots=s_n=1$, Formula (\ref{preciseweightedformula}) 
specializes to the identity 
\begin{eqnarray*}
\sum_{T\in\mathcal{HT}_k(n)}\prod_{e\in\mathcal E(T)}
(\sharp(e)-2)!=n^{k-1}S_1(n-1,k)(-1)^{n+k+1}
\end{eqnarray*}
which is analogous
to a Theorem of Husimi (see \cite{H} or \cite{GK}) expressing 
the total number
$$n^{k-1}S_2(n-1,k)$$
of elements in the set $\mathcal{HT}_k(n)$ of labelled hypertrees with $k$
hyperedges and $n$ vertices in terms of Stirling numbers of the second kind.

All these results are of course 
generalizations and variations of Cayley's theorem corresponding to 
the case $k=n-1$ and showing that there are 
$n^{n-2}$ labelled trees on $n$ vertices.

Observe that all these identities can also be deduced for example from 
Exercice 5.30 in \cite{St} using a well-known map between hypergraphs
and ordinary bipartite graphs.
\end{rem}


\section{Computational aspects and examples}\label{sectcomputexples}

The computation of $U(x,\sigma_1,\sigma_2,\dots)$ up to $o(x^n)$ is
straightforward using the recurrence relation
(\ref{eqtriangle}). For a given fixed numerical 
value of $\sigma_1$, the following
trick reduces memory requirement and speeds the computation up:
Setting
$$c_n(\sigma_1)=(\gamma_{n+1}(\sigma_1,n),\gamma_{n+2}(\sigma_1,n),
\dots,\gamma_{2n}(\sigma_1,n))$$ 
with $\gamma_i(\sigma_1,n)=\sum_{j=0}^{2n-i} t_{i,j}(n)(-\sigma_1)^j$
we have
$$U(x,\sigma_1,\sigma_2,\dots)=1-\sum_{n=1}^\infty
\langle c_n(\sigma_1),(\sigma_{n+1},\dots,\sigma_{2n})\rangle x^n$$
where $\langle a,b\rangle=\sum_{i\in I} a_ib_i$ for two finite-dimensional
vectors $a,b$ with coefficients indexed by a common 
finite set $I$.
The coefficients $\gamma_i(\sigma_1,n)$ of $c_n(\sigma_1)$ can be computed 
from the coefficients of $c_{n-1}(\sigma_1)$ by the formula
\begin{eqnarray}\label{formulagammai}
\gamma_i(\sigma_1,n)=
(i-2-\sigma_1)\gamma_{i-1}(\sigma_1,n-1)+(i-3)\gamma_{i-2}(\sigma_1,n-1)
\end{eqnarray}
with missing coefficients omitted in the case of $i=n+1$ or $i=2n$.

The coefficients of the first vectors $c_1(0),c_2(0),c_3(0),\dots$ are given 
by the rows of 
$$\begin{array}{rrrrr}
1\\1&1\\2&5&3\\
6&26&35&15\\
24&154&340&315&105\ ,\end{array}$$
see A112486 of \cite{OEIS}.


\subsection{The examples $U(x,-1,-1,-1,\dots)$ and $U(x,0,-1,-1,-1,\dots)$}
\label{sectexple111} 

The series
$$U(x,-1,-1,-1,-1,\dots)-1$$ 
is the generating series of the sequence 
$$S(n)=\sum_{i,j}t_{i,j}(n)$$
enumerating the sums of the triangles $T(n)$ defined by 
the integers $t_{i,j}(n)$. We have
\begin{eqnarray*}
&&(1+x)U(x,-1,-1,-1,-1,\dots)\\
&=&U(x,0,-2,-2,-2,-2,\dots)\\
&=&2U(x,0,-1,-1,-1,-1,\dots)-1
\end{eqnarray*}
where $U(x,0,-1,-1,-1,\dots)-1$ corresponds to the generating series
of the sequence 
$$s(n)=\sum_{i=n+1}^{2n} t_{i,0}(n)$$ 
starting as 
$$1, 2, 10, 82, 938, 13778, 247210, 5240338, 128149802, 3551246162,\dots,$$
cf. A112487 of \cite{OEIS},
and obtained by summing the integers
of the first column of the triangles $T(1),T(2),\dots$.
In particular, we have $2s(n)=S(n-1)+S(n)$ or equivalently
$$2\sum_{i=n+1}^{2n}t_{i,0}(n)=\sum_{i=n+1}^{2n}\sum_{j=0}^{2n-i}t_{i,j}(n)+
\sum_{i=n}^{2n-2}\sum_{j=0}^{2n-2-i}t_{i,j}(n-1)$$
for all $n\geq 2$.

\subsection{Examples satisfying differential equations}\label{sectdiffeq}

The recursive definition of the integers $t_{i,j}(n)$ implies easily that
specializations of the form
$$\sigma_n=c\frac{\prod_{k=1}^A(n+a_k)!}{\prod_{l=1}^B(n+b_l)!}z^{\lambda n+r},\ n\geq 2,$$
or 
$$\sigma_n=c\frac{\prod_{k=1}^A(n+a_k)!}{\prod_{l=1}^B(n+b_l)!}e^{(\lambda 
n+r)z},\ n\geq 2,$$
(with $\lambda\not=0$ and $b_i\not\in\{-2,-3,-4,-5,\dots\}$)
lead to differential equations with respect to $z$ for rational expressions
of $U(x,-y,\sigma_2,\sigma_3,\sigma_4,\dots)$. 
Such a series $U$ is analytic if $B>A$. 
We illustrate this with the following examples.
\subsubsection{$U(x,-y,-z^{2+r},-z^{3+r},-z^{4+r},\dots)$}
Setting $\sigma_n=-z^{n+r}$ for $n=2,3,\dots$ the series 
$f(z)=U(x,-y,\sigma_2,\sigma_3,\dots)-1$
satisfies formally the differential equation
$$f=xz\left(z^{r+1}+(y-(1+z)(1+r))f+z(1+z)\frac{df}{dz}\right)\ .$$
\subsubsection{$U\left(x,-y,-\frac{z^{2+r}}{(2+b)!},-\frac{z^{3+r}}{(3+b)!},-
\frac{z^{4+r}}{(4+b)!},\dots\right)$}
Setting $\sigma_n=-\frac{z^{n+r}}{(n+b)!}$ for $n=2,3,\dots$ the series 
$f(z)=U(x,-y,\sigma_2,\sigma_3,\dots)-1$
satisfies the differential equation
\begin{eqnarray*}
&&(b-r)(b-1-r)f+2(b-r)zf'+z^2f''\\
&=&\frac{xz^{2+r}}{b!}+xz(r-b-z(r+1)+(r-b)(r-y))f\\
&&+xz^2(b+y+z-2r)f'+xz^3f''
\end{eqnarray*}
\subsubsection{$U\left(x,-y,-e^{(2\lambda+r) z},-e^{(3\lambda+r) z},-e^{(4\lambda+r) z},\dots\right)$}
Setting $\sigma_n=-e^{(n\lambda+r)z}$ for $n=2,3,\dots$, the series
$$f(z)=U\left(x,-y,\sigma_2,\sigma_3,\sigma_4,\dots\right)-1$$
satisfies formally the differential equation
$$f=xe^{\lambda z}\left(e^{(\lambda+r)z}+\left(y-\left(1+\frac{r}{\lambda}
\right)\left(1+e^{\lambda z}\right)\right)f+\frac{1+e^{\lambda z}}{\lambda}
\frac{df}{dz}\right)\ .$$
 
\subsubsection{$U\left(x,-y,-(2+a)!e^{(2\lambda+r) z},-(3+a)!e^{(3\lambda+r) z},-(4+a!)e^{(4\lambda+r) z},\dots\right)$}
Setting $\sigma_n=-(n+a)!e^{(n\lambda+r)z}$ for $n=2,3,\dots$, the series
$f(z)=U\left(x,-y,\sigma_2,\sigma_3,\dots\right)-1$
satisfies formally the differential equation
\begin{eqnarray*}
f&=&\frac{xe^{\lambda z}}{\lambda^3}(\lambda(1+a)-r)\left(\lambda(\lambda(y-1)-r)+(\lambda+r)(r-\lambda(2+a))e^{\lambda z}\right)f\\
&&+\frac{xe^{\lambda z}}{\lambda^3}\left(\lambda(\lambda(y+a)-2r)+(3r^2-4\lambda r(a+1)+\lambda^2(a^2+a-1))e^{\lambda z}\right)f'\\
&&+\frac{xe^{\lambda z}}{\lambda^3}\left(\lambda+(2\lambda(1+a)-3r)e^{\lambda z}\right)f''+\frac{xe^{2\lambda z}}{\lambda^3}f'''+(2+a)!xe^{(2\lambda+r)z}\\
\end{eqnarray*}

\subsubsection{$U\left(x,-y,-\frac{e^{(2\lambda +r)z}}{(2+b)!},-\frac{e^{(3\lambda+r) z}}{(3+b)!},-\frac{e^{(4\lambda+r) z}}{(4+b)!}\dots\right)$}
Setting $\sigma_n=-\frac{e^{(n\lambda +r)z}}{(n+b)!}$ for $n=2,3,\dots$
the series $f(z)=U\left(x,-y,\sigma_2,\sigma_3,\dots\right)-1$
satisfies the differential equation
\begin{eqnarray*}
&&(\lambda b-r)(\lambda(b-1)-r)f+
(\lambda(2b-1)-2r)f'+f''\\
&=&x\frac{\lambda^2}{b!}e^{(2\lambda+r)z}+
xe^{\lambda z}\left((\lambda b-r)(\lambda(y-1)-r)-\lambda(\lambda+r)e^{\lambda z}\right)
f\\
&&+xe^{\lambda z}\left(\lambda(b-1+y)-2r+\lambda e^{\lambda z}\right) f'
+xe^{\lambda z}f''
\end{eqnarray*}

\begin{rem}
The recursion relation (\ref{eqtriangle})
gives rise to partial differential equations for generating
series of $t_{i,j}(n)$ which are exponential 
with respect to $j$ and/or $n$.
\end{rem}


\subsection{A family of rational examples}
\begin{prop}\label{proprationalexple} Let $\sigma_1,\sigma_2,\dots$ be a sequence of complex numbers 
of the form $\sigma_n=(-1)^nP(n)$ for all $n\geq A$ where $A$ is some
natural integer and where $P(s)\in\mathbb C[s]$ is a polynomial.
Then $U(x,\sigma_1,\sigma_2,\dots)$ is a rational series.
\end{prop}

{\bf Proof} Let $d$ denote the degree of $P$. Applying identity 
(\ref{functionaleqU}) of 
Proposition~\ref{propfunctionalequationexistence} iteratively
$d+1$ times we get a series of the form $U(x,\tilde \sigma_1,\tilde \sigma_2,
\dots,\tilde \sigma_{A+d+2},0,0,0,\dots)$ which is a polynomial.\hfill
$\Box$

As an illustration we consider the series 
$U(x,y,1,-1,1,\dots)$. 
Proposition~\ref{propfunctionalequationexistence} shows 
$$(1-xy)U(x,y,1,-1,1,-1,\dots)=U(x,1+y,1+y,0,0,\dots)=1-(1+y)x\ .$$
We have thus $U(x,y,1,-1,1,\dots)=1-\frac{x}{1-xy}$.

\subsection{Coefficients of $U(x,\sigma_1,P(2),P(3),P(4),\dots)$}

\begin{prop} Let $P(s)\in\mathbb C[s]$ be a polynomial of degree $d$.
There exist constants $\alpha_0,\dots,\alpha_d\in\mathbb C$ such that
$$[x^n]U(x,\sigma_1,P(2),P(3),P(4),\dots)=
\sum_{h=0}^d\alpha_h[x^{n+h}]U(x,\sigma_1,1,1,1,1,1,\dots)$$
for all $n\geq 1$ with 
$[x^n]U$ denoting the coefficient of $x^n$ in 
the series $U$.
\end{prop}

{\bf Proof} The proof is by induction on $d$ and holds certainly for $d=0$.
Setting $\gamma_i(n)=\sum_{j=n+1}^{2n}t_{i,j}(n)(-\sigma_1)^j$, formula
(\ref{formulagammai}) implies
\begin{eqnarray*}
0&=&-i^d\gamma_i(n+1)+i^d(i-2-\sigma_1)\gamma_{i-1}(n)+i^d(i-3)\gamma_{i-2}
(n)\\
&=&-i^d\gamma_i(n+1)+(i-1)^{d+1}\gamma_{i-1}(n)+(i-2)^{d+1}\gamma_{i-2}(n)+\\
&&\ +Q_1(i-1)\gamma_{i-1}(n)+Q_2(i-2)\gamma_{i-2}(n)
\end{eqnarray*}
where $Q_1$ and $Q_2$ are polynomials of degree $\leq d$. Fixing $n$ 
and summing over $i$ we get 
\begin{eqnarray}\label{sumoveri}
2\sum_{i=n+1}^{2n}i^{d+1}\gamma_i(n)=\sum_{i=n+2}^{2n+2}i^d\gamma_i(n+1)
-\sum_{i=n+1}^{2n}(Q_1+Q_2)(i)\gamma_i(n)\ .
\end{eqnarray}
The right side of (\ref{sumoveri}) equals now
$$[x^{n+1}]U(x,\sigma_1,2^d,3^d,4^d,\dots)
-[x^n]U(x,\sigma_1,(Q_1+Q_2)(2),
(Q_1+Q_2)(3),\dots)\ .$$
It is thus by induction on $d$ a linear combination of the coefficients 
of $x^n,\dots,x^{n+d+1}$ in $U(x,\sigma_1,1,1,1,\dots)$.
This proves the result for $U(x,\sigma_1,2^{d+1},3^{d+1},\dots)$. 
The general induction step follows by remarking that all coefficients 
of strictly positive degree in $x$ of $U(x,\sigma_1,\sigma_2,\dots)$
are linear in $\sigma_2,\sigma_3,\dots$.\hfill$\Box$


\section{Conjectural asymptotics for $s(1),s(2),\dots$}
\label{sectasympt}

Computations with a few thousand values of $s(n)$ suggest the following
asymptotic formula for 
the integral sequence $s(n)=\sum_{i=n+1}^{2n}t_{i,0}(n)$:

\begin{conj}\label{conjasympt}
There exists a sequence 
$A_0,A_1,\dots$ of rational polynomials 
$A_i(x)\in\mathbb Q[x]$ with $A_i$ of degree $i$ such that
$$s(n)=\frac{n^{n-1}}{(1-\log 2)^{n-1/2}e^n}
\left(\sum_{k=0}^m\frac{A_k(1-\log 2)}{n^k}+o(n^{-m})\right)$$
for all $m\in\mathbb N$.
\end{conj}

The first few polynomials $A_0,A_1,\dots$ are
\begin{eqnarray*}
A_0&=&1\\
A_1&=&\frac{11}{24}-\frac{x}{12}\\
A_2&=&\frac{265}{1152}-\frac{47x}{288}+\frac{x^2}{288}\\
A_3&=&\frac{48703}{414720}-\frac{3649x}{13824}+\frac{107x^2}{6912}
+\frac{139x^3}{51840}\\
A_4&=&\frac{2333717}{39813120}-\frac{2019163x}{4976640}
+\frac{16489x^2}{331776}+\frac{26549x^3}{1244160}-\frac{571x^4}{2488320}\\
A_5&=&\frac{38180761}{1337720832}-\frac{293093189x}{477757440}
+\frac{16859263x^2}{119439360}+\\
&&\ +\frac{6752203x^3}{59719680}-\frac{170729x^4}{59719680}-
\frac{163879x^5}{209018880}\\
\end{eqnarray*}
The coefficients $B_k$ of the formal power series
$\sum_{k=1}^\infty B_k(x)t^k=\log\left(\sum_{k=0}^\infty
A_k(x)t^k\right)$ seem to be simpler and start as
\begin{eqnarray*}
B_1&=&\frac{11}{24}-\frac{x}{12}\\
B_2&=&\frac{1}{8}-\frac{x}{8}\\
B_3&=&\frac{127}{2880}-\frac{3x}{16}+\frac{x^2}{288}
+\frac{x^3}{360}\\
B_4&=&\frac{1}{64}-\frac{9x}{32}+\frac{11x^2}{576}+\frac{x^3}{48}\\
B_5&=&\frac {221}{40320}-\frac{27x}{64}+\frac{41x^2}{576}+
\frac{1381x^3}{12960}-\frac {x^4}{1440}-\frac{x^5}{1260}
\end{eqnarray*}

\begin{rem}
The constant $1-\log 2=.30685281944\dots$ appearing in 
Conjecture \ref{conjasympt} seems also to 
be related to the index $m_n$ such that $t_{m_n,0}(n)=
\max_i(t_{i,0}(n))$ with $m_n$ given asymptotically by 
$\frac{n}{2(1-\log 2)}$. Moreover, we have seemingly
$\lim_{n\rightarrow\infty}\frac{t_{m_n,0}(n)\ \sqrt{n}}{s(n)}\sim
.87$ (and the numbers $t_{i,0}(n)$, suitably rescaled, should 
satisfy a central limit Theorem).
\end{rem}


\section{Modular properties of the sequence $s(1),s(2),\dots$}
\label{sectmodular}

\begin{prop}\label{proprationalmodp}
The series $U(x,\sigma_1,\sigma_2,\dots)\in\mathbb F_p[[x]]$
is rational if $\sigma_1,\sigma_2,\dots$ is an ultimately periodic 
sequence of elements in $\mathbb F_p$.
\end{prop}

\noindent{\bf Proof} Up to addition of a polynomial to
$U=U(x,\sigma_1,\sigma_2,\dots)$ we can suppose that
$\sigma_2,\sigma_3,\sigma_4,\dots$ is periodic with period $k$.
We set $\tilde \sigma_i=\sigma_i$ for $i\geq 2$ and extend $\tilde\sigma_2,
\tilde\sigma_3,\dots,$ to a $k-$periodic sequence indexed by $\mathbb Z$.
We suppose first $\sigma_1\not=0$ in $\mathbb F_p$. Using the 
identity $(-\sigma_1)^{p-1}=1$ and periodicity of the sequence
$(\tilde\sigma_i)_{i\in \mathbb Z}$, we have
\begin{eqnarray*}
U&=&1-\sum_{n=1}^\infty x^n\sum_{i=0}^{pk-1}\sum_{\alpha\in\mathbb Z}
\tilde\sigma_{i+\alpha kp}
\sum_{j=0}^{p-2}\sum_{\beta\in\mathbb Z}t_{i+\alpha kp,j+\beta(p-1)}(n)(-\sigma_1)^{j+
\beta(p-1)}\\
&=&1-\sum_{n=1}^\infty x^n\sum_{i=0}^{pk-1}\tilde\sigma_{i}
\sum_{j=0}^{p-2}(-\sigma_1)^j\left(
\sum_{\alpha\in\mathbb Z}\sum_{\beta\in\mathbb Z}t_{i+\alpha kp,j+\beta(p-1)}(n)\right)\ .
\\
\end{eqnarray*}
Since the recurrence relations (\ref{eqtriangle}) define the elements
$t_{i,j}(2),t_{i,j}(3),\dots$ correctly for arbitrary indices $i,j\in\mathbb Z$
we have
\begin{eqnarray*}
&&\sum_{\alpha\in\mathbb Z}\sum_{\beta\in\mathbb Z}t_{i+\alpha kp,j+\beta(p-1)}(n)\\
&=&(i-2)
\sum_{\alpha\in\mathbb Z}\sum_{\beta\in\mathbb Z}t_{i-1+\alpha kp,j+\beta(p-1)}(n-1)\\
&&\ +\sum_{\alpha\in\mathbb Z}\sum_{\beta\in\mathbb Z}t_{i-1+\alpha kp,j-1+\beta(p-1)}(n-1)\\
&&\ +(i-3)
\sum_{\alpha\in\mathbb Z}\sum_{\beta\in\mathbb Z}t_{i-2+\alpha kp,j+\beta(p-1)}(n-1)
\end{eqnarray*}
for $n\geq 2$. Setting 
$$\tilde t_{i,j}(n)\equiv
\sum_{\alpha\in\mathbb Z}\sum_{\beta\in\mathbb Z}t_{i+\alpha kp,j+\beta(p-1)}(n-1)$$
for $0\leq i<kp$ and $0\leq j<p-1$, the elements
$\tilde t_{i,j}(n)$ of $\mathbb F_p$  satisfy the recursion 
relation (\ref{eqtriangle}) with indices considered modulo $kp$ 
for $i$ and modulo $p-1$ for $j$. Since the $kp(p-1)$ elements 
$\tilde t_{i,j}(n)$ of the finite field $\mathbb F_p$ 
depend affinely on the $kp(p-1)$ elements
$\tilde t_{i,j}(n-1)$ for $n\geq 2$, finiteness of the set $(i,j)$ of
indices implies the existence of an integer $l$ such that 
$\tilde t_{i,j}(n+l)=t_{i,j}(n)$ for all
sufficiently large $n$ and for all possible indices $i$ and $j$. 
This implies easily that the coefficients of $U$ are ultimately
periodic and ends the proof for $\sigma_1\not=0$.

The case $\sigma_1=0$ involves only the integers 
$t_{i,0}(n)$ and their analogues $\tilde t_{i,0}(n)$
with indices in the finite set $\{0,\dots,pk-1\}$. 
Details are similar to the previous case and left to the reader.
\hfill$\Box$

The first non-trivial case of Proposition \ref{proprationalmodp}
is perhaps given by the generating series $U(x,0,-1,-1,-1,\dots)$
with coefficients of $U(x,0,-1,-1,\dots)-1$ given by the sequence 
$$s(n)=\sum_{i=n+1}^{2n}t_{i,0}(n)$$
obtained by summing all coefficients in the
first column of the triangular arrays
$T(1),T(2),\dots$. 

\begin{conj} There exists a sequence 
\begin{eqnarray*}
&&\alpha_0=-1,\alpha_1=2,\alpha_2=0,\alpha_3=\frac{1}{3},\alpha_4={\frac {5}{18}},\alpha_5={\frac {149}{540}},\alpha_6={\frac {553}{2025}},\\
&&\alpha_7={\frac {1849741}{6804000}},\alpha_8={\frac {775167119}{2857680000}},\alpha_9={\frac {
325214957371}{1200225600000}},\dots
\end{eqnarray*}
of rational numbers such that 
$$\left(1+x^{p-1}\right)\sum_{n=1}^\infty s(n)x^n\equiv x+\sum_{n=0}^{p-2}
\alpha_nx^{p-n}\pmod p$$
for every prime number $p$.
\end{conj}


\begin{conj}
The rational sequence $\alpha_0,\alpha_1,\dots$ has
an asymptotic expansion given by 
$$\alpha_n\sim
\sum_{k=1}^\infty \frac{k^{k-n}}{k!}\left(\frac{2}{e^2}\right)^k$$
and converges with limit given by
$2e^{-2}=.27067056647322538378799\dots$.

The error term 
$$\epsilon_n
=\alpha_n-\sum_{k=1}^\infty \frac{k^{k-n}}{k!}\left(\frac{2}{e^2}\right)^k
$$
is given by
$$\epsilon_n= \frac{(-1)^{n+1}}{(1-\log 2)s(n+1)}\left(\sum_{k=0}^m
\frac{\gamma_{2k}(1-\log 2)}{n^{2k}}+o\left(n^{-2m-1}\right)\right)$$
where $\gamma_{2k}(x)\in \mathbb Q[x]$ is a polynomial of degree 
at most $2k$. The first few polynomials are given by
\begin{eqnarray*}
\gamma_0&=&1\\
\gamma_2&=&-\frac{x}{12}\\
\gamma_4&=&-\frac{x}{48}+\frac{x^2}{48}+\frac{x^3}{40}\\
\gamma_6&=&-\frac{x}{192}+\frac{5x^2}{96}+\frac{193x^3}{864}-\frac{x^4}{72}
-\frac{5x^5}{252}\\
\end{eqnarray*}
\end{conj}

\section{Integer sequences obtained as weighted sums of the numbers $t_{i,j}(n)$}\label{sectinttijn}

The sequence $q_i(x,y)$ defined by
$$q_i(x,y)=
\sum_{n,j\geq 0} t_{i,j}(n)x^ny^j=
\sum_{n=\lceil i/2\rceil}^{i-1}x^n\sum_{j=0}^{n-1}t_{i,j}(n)y^j$$
is given by $q_1(x,y)=0,\ q_2(x,y)=x$ and by the recursion relation
\begin{eqnarray*}
q_i(x,y)=
x\left((i-2+y)q_{i-1}(x,y)+(i-3)q_{i-2}(x,y)\right)
\end{eqnarray*}
for $i\geq 3$.
The following table lists the first few non-zero coefficients $q_2(x,y),
q_3(x,y),\dots$ (up to normalizations) 
and the seemingly corresponding sequences of \cite{OEIS}
(which have often interesting combinatorial interpretations) 
for a few specializations:
$$\begin{array}{|l|c|c|}
\hline
q_i(-3,-2/3)(-1)^{i+1}/3&1,1,1,1,1,1,\dots&A12\\
\hline
q_i(-3,1/3)(-1)^{i+1}/3&1,4,25,226,2713,40696,\dots&A10845\\
\hline
q_i(-2,-3/2)(-1)^{i+1}/2&1,-1,-3,-5,-7,-9,-11,\dots&\\
\hline
q_i(-2,-1/2)(-1)^{i+1}/2&1,1,1,1,1,1,\dots&A12\\
\hline
q_i(-2,1/2)(-1)^{i+1}/2&1,3,13,79,633,6331,\dots&A10844\\
\hline
q_i(-1,-1)(-1)^{i+1}&1,0,-1,-2,-3,-4,-5,\dots&\\
\hline
q_i(-1,0)(-1)^{i+1}&1,1,1,1,1,1,\dots&A12\\
\hline
q_i(-1,1)(-1)^{i+1}&1,2,5,16,65,326,1957,\dots&A522\\
\hline
q_i(-1,2)(-1)^{i+1}&1,3,11,49,261,1631,\dots&A1339\\
\hline
q_i(-1,3)(-1)^{i+1}&1,4,19,106,685,5056,\dots&A82030\\
\hline
q_i(1,-2)(-1)^i&1,1,1,1,1,1,\dots&A12\\
\hline
q_i(1,-1)&1,0,1,2,9,44,265,1854,\dots  &A166\\
\hline
q_i(1,0)&1,1,3,11,53,309,2119,
\dots&A255\\
\hline
q_i(1,1)&1,2,7,32,181,1214,\dots&A153\\
\hline
q_i(1,2)&1,3,13,71,465,3539,\dots&A261\\
\hline
q_i(1,3)&1,4,21,134,1001,8544,\dots&A1909\\
\hline
q_i(2,-7/2)(-1)^i/2&1,5,17,37,65,101,145,197,\dots&A53755\\
\hline
q_i(2,-5/2)(-1)^i/2&1,3,5,7,9,11,\dots&A5408\\
\hline
q_i(2,-3/2)(-1)^i/2&1,1,1,1,1,1,\dots&A12\\
\hline
q_i(2,-1)/2&1,0,2,8,60,544,6040,\dots&A53871\\
\hline
q_i(2,-1/2)/2&1,1,5,29,233,2329,\dots&A354\\
\hline
q_i(3,-4/3)(-1)^i/3&1,1,1,1,1,1,\dots&A12\\
\hline
q_i(3,-1)/3&1,0,3,18,189,2484,\dots&A33030\\
\hline
\end{array}$$

The sequences A255, A153, A261 and A1909 corresponding to $q_i(1,0),q_i(1,1),
q_i(1,2)$ and $q_i(1,3)$ can seemingly also be obtained by 
considering the weighted sums
$$s_i=\sum_{j,n\geq 0}{n-j+k-1\choose k}t_{i,j}(n)(-1)^j$$
for $k=1,2,3$ and $4$.

The above table contains a few instances of the identities
$$1=q_i\left(\kappa,-\frac{\kappa+1}{\kappa}\right)\frac{(-1)^i}{\kappa}$$
and
$$1+(i-2)\kappa=q_i\left(\kappa,-\frac{2\kappa+1}{\kappa}\right)
\frac{(-1)^i}{\kappa}$$
which hold for $\kappa\not=0$ and for $i=2,3,4,\dots$ and which 
can easily be proven by induction. (These two examples generalize
probably to $q_i\left(\kappa,-\frac{\lambda\kappa+1}{\kappa}\right)\frac{(-1)^i}{\kappa}=P_\lambda(i,\kappa),\ \lambda=1,2,3,\dots,\ \kappa\not=0,\ 
i=2,3,4,\dots$ for $P_\lambda$ a suitable polynomial function of $\kappa$ 
and $i$.)

Another identity is given by the family of weighted examples 
$$1=(-1)^ik!(i-1+k)\sum_{j,n\geq 0}t_{i,j}(n)\frac{(-1)^j}{(n+k-j)!}$$
for all $k\in \mathbb N$ and for all $i\geq 2$.

A few other interesting weighted examples (there are probably many more) 
are given by the following formulae
\begin{eqnarray*}
a_i&=&\sum_{j,n\geq 0}(n-j)t_{i,j}(n)(-2)^j\ ,\\
b_i&=&\frac{1}{4}\sum_{j,n\geq 0}(n-1-j)t_{i,j}(n)2^n\left(\frac{-3}{2}\right)^j\ ,\\
c_i&=&\sum_{j,n\geq 0}t_{i,j}(n)(n-1-j)!(-1)^j\ ,\\
d_i&=&(-1)^{i+1}\sum_{j,n\geq 0}t_{i,j}(n)(n-1-j)!(-1)^n\ ,\\
e_i&=&\sum_{j,n\geq 0}t_{i,j}(n)\frac{(-1)^j}{(n-1-j)!}\ ,\\
f_i&=&(i-1)!\sum_{n\geq 1}\frac{t_{i,0}(n)}{(n-1)!}\ .
\end{eqnarray*}
Their initial coefficients (with leading zeros omitted)
and the seemingly corresponding sequences of \cite{OEIS} are as
follows:
$$\begin{array}{|c|l|c|}
\hline
a_i&1,0,0,1,1,8,36,229,1625,\dots&A757\\
\hline
b_i&1,0,5,24,209,2120,\dots&A120765\\
\hline
c_i&1,0,3,26,453,11844,\dots&A89041\\
\hline
d_i&1,2,7,52,749,17686,\dots&A46662\\
\hline
e_i&1,0,0,0,0,0,0,\dots&A7\\
\hline
f_i&1,2,12,84,820,9540,\dots&A179495\\
\hline
\end{array}$$
Most of the proofs are probably easy: The corresponding sections of
\cite{OEIS} give information concerning generating functions which 
can be applied to differential equations 
analogous to those of Section \ref{sectdiffeq}.

\section{Coverings}\label{sectcov}

Coverings and packings are dual notions. We discuss here a few 
aspects of the theory of coverings in relation with packings 
by generic families.

A {\em (left-)covering} with parts $\mathcal S_1,\dots,\mathcal S_n$
of a group $G$ is a vector $(g_1,\dots,g_n)$ such that 
$G=\cup_{j=1}^n g_j\mathcal S_j$.

A covering of a finite group $G$ with $N$ elements by 
non-empty subsets $\mathcal S_1,\dots,\mathcal S_n$ exists of 
course always if $n\geq N$.

We are interested in large collections of subsets 
$\mathcal S_1,\dots,\mathcal S_n$ in a finite group $G$ of order $N$ 
such that the sets 
$\mathcal S_1,\dots,\mathcal S_n$ (or more precisely, suitable translates)
cover $G$ and the number 
of all coverings depends only on the cardinalities of $\mathcal S_1,
\dots,\mathcal S_n$ (and of $N$) for every family $\mathcal S_1,\dots,
\mathcal S_n$ in the collection. 

Three such collections can be described as follows:

Start with a family $\mathcal S_1,\dots,S_n$ which is generic for packings 
and add $N-\sum_{j=1}^n\sharp(\mathcal S_j)$ singletons.
Coverings of $G$ by such families are \lq\lq tight\rq\rq and
essentially in one-to-one correspondence (except for a factor 
$\left(N-\sum_{j=1}^n\sharp(\mathcal S_j)\right)!$ accounting for all 
permutations of the added singletons) with packings by $\mathcal S_1,\dots,
\mathcal S_n$.

The second family is obtained by adding 
$N+n-1-\sum_{j=1}^n\sharp(\mathcal S_j)$ singletons to a 
family $\mathcal S_1,\dots,S_n$ which is generic for packings.
The fact that the number of associated coverings depends only on
all involved cardinalities is similar to the proof of Proposition 
\ref{propexistence} given in Section \ref{sectexistence}. 
The proof needs probably computations with the full M\"obius function. 
I do not know if there is an efficient way for computing the number
of associated coverings or if there is 
a nice formula similar to the one assciated to enumerations of packings. 

There is a further variation 
on this theme: Given an arbitrary natural integer $a$ one can 
consider adding 
$N+a-\sum_{j=1}^n\sharp(\mathcal S_j)$ singletons to 
a family $\mathcal S_1,\dots,S_n$ which is generic for packings.
For every natural integer $a$, such a family has 
the property that the number of associated coverings 
depends only on all involved cardinalities. The choices $a=0$ corresponding 
to the first family and $a=n-1$ corresponding to the second family 
are however natural in this context. Indeed, since the intersection graph 
of the sets $g_1\mathcal S_1,\dots,g_n\mathcal S_n$ is a hyperforest,
the union $\cup_{j=1}^ng_j\mathcal S_j$ contains at least $\sum_{j=1}^n\sharp
(\mathcal S_j)-(n-1)$ elements. 
This leaves at most $N+n-1-\sum_{j=1}^n\sharp
(\mathcal S_j)$ missing elements which can be covered using the additional
singletons.

A third rather trivial family is given by 
considering complements $G\setminus\mathcal S_1,\dots,G\setminus\mathcal S_n$
where $\mathcal S_1,\dots,\mathcal S_n$ is a generic family for packings 
in $G$ having at least two parts. The number of coverings of such a family
is easy to compute and given by $N^n-N\prod_{j=1}^n\sharp(\mathcal S_j)$.

It would perhaps be interesting to have other (and hopefully more exotic)
families of examples.

\noindent{\bf Acknowledgements.} 
I thank Pierre de la Harpe for helpful comments
and two anonymous referees for their careful work and useful remarks.

\noindent Roland BACHER, Universit\'e Grenoble I, CNRS UMR 5582, Institut 
Fourier, 100 rue des maths, BP 74, F-38402 St. Martin d'H\`eres, France.

\noindent e-mail: Roland.Bacher@ujf-grenoble.fr


\begin{thebibliography}{99}

\bibitem{GK} I.M. Gessel, L.H. Kalikow, {\it
Hypergraphs and a functional equation of Bouwkamp and de Bruijn.}
J. Combin. Theory Ser. A 110 (2005), no. 2, 275--289. 

\bibitem{H} K. Husimi, {\it Note on Mayer's theory of cluster integrals},
Journal of Chemical Physics {\bf 18} (1950), 682--684.

\bibitem{MM} J. McCammond, J. Meier, {\it
The hypertree poset and the $l^2$-Betti numbers of the motion group of the trivial link.} Math. Ann. {\bf 328} (2004), no. 4, 633--652. 

\bibitem{OEIS} N.J.A. Sloane, The On-Line Encyclopedia of Integer Sequences, published electronically at http://oeis.org.

\bibitem{St} R. P. Stanley, Enumerative Combinatorics, Volume I, 
Cambridge University Press (1997).

\end{thebibliography}
\end{document}